\documentclass[12pt]{amsart}
\usepackage[]{amscd,amsfonts,amsmath,amsxtra,amssymb}
\usepackage{graphics}
\usepackage[all]{xy}
\usepackage{hyperref}
\usepackage{eucal}
\usepackage{rotating}
\usepackage[T1]{fontenc}
\newtheorem{sub}{}[section]
\newtheorem{subsub}{}[sub]
\usepackage{rotating}
\def\ov#1{\overline{#1}}

\def\coker{\mathop{\rm coker}\nolimits}
\def\Hom{\mathop{\rm Hom}\nolimits}

\def\HHom{\mathop{\mathcal Hom}\nolimits}

\def\Ext{\mathop{\rm Ext}\nolimits}
\def\EExt{\mathop{\mathcal Ext}\nolimits}

\def\Pic{\mathop{\rm Pic}\nolimits}

\def\Aut{\mathop{\rm Aut}\nolimits}

\def\End{\mathop{\rm End}\nolimits}

\def\EEnd{\mathop{\mathcal End}\nolimits}

\def\imm{\mathop{\rm im}\nolimits}
\def\deg{\mathop{\rm deg}\nolimits}
\def\rg{\mathop{\rm rg}\nolimits}
\def\rk{\mathop{\rm rk}\nolimits}
\def\Deg{\mathop{\rm Deg}\nolimits}

\def\spec{\mathop{\rm spec}\nolimits}

\def\lra{\longrightarrow}

\def\sigg{\mathop{\hbox{$\displaystyle\sum$}}\limits}

\def\hfl#1#2{\smash{\mathop{\ \hbox to 12mm{\rightarrowfill}}
\limits^{\scriptstyle#1}_{\scriptstyle#2} \ }}
\def\hflb#1#2{\smash{\mathop{\hbox to 12mm{\leftarrowfill}}
\limits^{\scriptstyle#1}_{\scriptstyle#2}}}

\def\m#1{{\hbox{$#1$}}}

\def\ot{\otimes}
\def\og{\leavevmode\raise.3ex\hbox{$\scriptscriptstyle\langle\!\langle$}}
\def\fg{\leavevmode\raise.3ex\hbox{$\scriptscriptstyle\,\rangle\!\rangle$}}
\def\span#1{\langle#1\rangle}

\def\nsp{\lbrace 0\rbrace}

\def\dsp{\displaystyle}
\def\Ssect#1#2{\pagebreak[3]\begin{sub}\label{#2}{\sc\small\small
#1}\rm\medskip}

\def\sepsec{\vskip 1.8cm}
\def\sepsub{\vskip 1cm}

\def\sepprop{\vskip 0.5cm}

\def\xmat#1{\[\xymatrix{#1}\]}
\def\flinc{\ar@{^{(}->}}
\def\fleq{\ar@{=}}
\def\flon{\ar@{->>}}
\def\fmaps{\ar@{|-{>}}}
\def\fflat{\ar@{-}}
\def\fimpl{\ar@{=>}}
\def\Nligne{\hfil\break}

\def\ED{\vskip 1cm\end{document}}

\newcommand{\M}{{\mathbb M}}

\newcommand{\C}{{\mathbb C}}

\renewcommand{\P}{{\mathbb P}}

\newcommand{\E}{{\mathbb E}}

\newcommand{\U}{{\mathbb U}}

\renewcommand{\L}{{\mathbb L}}

\newcommand{\kc}{{\mathcal C}}
\newcommand{\kd}{{\mathcal D}}
\newcommand{\ke}{{\mathcal E}}
\newcommand{\kf}{{\mathcal F}}

\newcommand{\ki}{{\mathcal I}}

\newcommand{\kl}{{\mathcal L}}

\newcommand{\kn}{{\mathcal N}}
\newcommand{\ko}{{\mathcal O}}

\newcommand{\ku}{{\mathcal U}}
\newcommand{\kv}{{\mathcal V}}

\newcommand{\kx}{{\mathcal X}}

\newcommand{\kz}{{\mathcal Z}}
\begin{document}

\def\refname{References}
\def\contentsname{Summary}
\def\proofname{Proof}
\def\abstractname{Resume}

\author{Jean--Marc Dr\'{e}zet}
\address{
Institut de Math\'ematiques de Jussieu - Paris Rive Gauche\\
Case 247\\
4 place Jussieu\\
F-75252 Paris, France}
\email{jean-marc.drezet@imj-prg.fr}
\title[{Non-reduced moduli spaces of sheaves}] {Non-reduced moduli spaces of 
sheaves on multiple curves}

\begin{abstract} Some coherent sheaves on projective varieties have a non 
reduced versal deformation space. For example, this is the case for most 
unstable rank 2 vector bundles on $\P_2$ (cf. \cite{st}). In particular, it may 
happen that some moduli spaces of stable sheaves are non reduced. 

We consider the case of some sheaves on ribbons (double structures on smooth 
projective curves): the quasi locally free sheaves of rigid type. Le $E$ be 
such a sheaf.

-- Let $\ke$ be a flat family of sheaves containing $E$. We find that it is a 
reduced deformation of $E$ when some canonical family associated to $\ke$ is 
also flat.

-- We consider a deformation of the ribbon to reduced projective curves with 
two components, and find that $E$ can be deformed in two distinct ways to 
sheaves on the reduced curves. In particular some components $\bf M$ of the 
moduli spaces of stable sheaves deform to two components of the moduli spaces 
of sheaves on the reduced curves, and $\bf M$ appears as the ``limit'' of 
varieties with two components, whence the non reduced structure of $\bf M$.
\end{abstract}

\maketitle
\tableofcontents


\section{Introduction}

Let $X$ be a projective variety over $\C$ and $E$ a coherent sheaf on $X$. Let 
\m{(V_E,v_0,\ke_E,\alpha)} be a semi-universal deformation of $E$ 
(\m{(V_E,v_0)} is a the germ of an analytic variety, \m{\ke_E} is a coherent 
sheaf on \m{X\times V_E} and \ \m{\alpha:\ke_{E,s_0}\to E} \ is an isomorphism, 
cf. \cite{si_tr}). It may happen that \m{V_E} is not reduced at \m{s_0}, for 
example in the case of unstable rank-2 vector bundles on \m{\P_2} (cf. 
\ref{P2}), or for some sheaves on ribbons. It is then natural to ask {\em why} 
\m{V_E} is not reduced. 

{\em Notation:} In this paper, an {\em algebraic variety} is a quasi-projective 
scheme over $\C$.
\sepsub

\Ssect{Non-reduced deformations of sheaves}{n-r-d}

\begin{subsub} The general case -- \rm
If $\kf$ is a flat family of sheaves parametrised by the germ \m{(T,t_0)} of an 
analytic variety, and if \m{\kf_{t_0}\simeq E}, there is a morphism \ \m{\phi: 
(T,t_0)\to(V_E,v_0)} \ such that \ \m{(I_X\times\phi)^*(\ke_E)\simeq\kf}, with 
a uniquely defined tangent map
\[T_{t_0}T\lra T_{v_0}V_E=\Ext^1_{\ko_X}(E,E) \ . \]
If $T$ is reduced then $\phi$ can be factorized to \m{V_{E,reg}}, the reduced 
germ associated to \m{V_E}:
\xmat{\phi:T\ar[r]^-{\phi'} & V_{E,reg}\flinc[r] & V_E \ . }
It is then natural to ask:
\begin{enumerate}
\item[(i)] What is the tangent space \ 
$T_{v_0}V_{E,reg}\subset\Ext^1_{\ko_X}(E,E)$ ?
\item[(ii)] Under which conditions on $\kf$, when $T$ is non reduced, is $\phi$ 
a morphism to $V_{E,reg}$ ?
\end{enumerate}
and the more vague question
\begin{enumerate}
\item[(iii)] Why can \m{V_E} be non reduced ? 
\end{enumerate}
The problem can also be stated in terms of moduli spaces of stable sheaves: let 
\m{\ko_X(1)} be an ample line bundle on $X$ and \m{P_E} the Hilbert polynomial 
of $E$. Let \m{\bf M} be the moduli space of sheaves on $X$, stable with 
respect to \m{\ko_X(1)}, and with Hilbert polynomial \m{P_E}. 
Suppose now that $E$ is stable. We can then ask
\begin{enumerate}
\item[(i')] What is the tangent space \ 
$T_E{\bf M}_{reg}\subset T_E{\bf M}$ ?
\item[(ii')] If $\kf$ is a flat family of stable sheaves of Hilbert polynomial 
$P_E$, parametrised by a variety $T$, we get a morphism $\phi:T\to{\bf M}$. 
Under which conditions on $\kf$, when $T$ is non reduced, is 
$\phi$ a morphism to ${\bf M}_{reg}$ ?
\item[(iii')] Why can $\bf M$ be non reduced ?
\end{enumerate}
\end{subsub}

For questions (ii), in the case of unstable rank-2 vector bundles on \m{\P_2} 
and of some sheaves on multiple curves studied here, it appears that the 
bundles (or sheaves) have a canonical extra structure, and that $\phi$ is a 
morphism to \m{V_{E,reg}} when this structure is also flat over \m{V_E}. The 
case unstable rank-2 vector bundles on \m{\P_2} is briefly recalled in \ref{P2}.

\sepprop

\begin{subsub} The case of sheaves on ribbons -- Good families of sheaves -- \rm
In this paper we will investigate the case of some sheaves on {\em ribbons} 
(i.e double structures on smooth projective curves). Let $Y$ be a ribbon and \ 
\m{C=Y_{reg}}. For vector bundles \m{V_E} is smooth (because in this case \ 
\m{\Ext^2_{\ko_Y}(E,E)=\nsp}). We will consider {\em quasi-locally free sheaves 
of rigid type} (defined and studied in \cite{dr4}), i.e. coherent sheaves 
locally isomorphic to \ \m{a\ko_Y\oplus\ko_C} \ for some integer $a$. 
Deformations of these sheaves are also quasi-locally free sheaves of rigid 
type, and \m{\deg(E_{|C})} is invariant by deformation.
\end{subsub}

For these sheaves, for (i) we have from \cite{dr4}
\[T_{v_0}V_{E,reg} \ = \ H^1(\End(E)) \ . \]
Let $L$ be the ideal sheaf of $C$ in $Y$. It is a line bundle on $C$. We will 
prove that
\[T_{v_0}V_E/T_{v_0}V_{E,reg} \ \simeq \ H^0(\EExt^1_{\ko_Y}(E,E)) \ \simeq \ 
H^0(L^*) \ . \]
For question (ii), let $Z$ be an algebraic variety and $\kf$ a coherent sheaf 
on \m{Z\times Y}, flat on $Z$. Let \m{z_0\in Z} be such that \m{\kf_{z_0}\simeq 
E}. So we get a morphism
\[\phi:(Z,z_0)\lra(V_E,v_0)\]
(where \m{(Z,z_0)} is the germ defined by $Z$ and \m{z_0}).
We say that $\ke$ is a {\em good family} if \m{\kf_{|Z\times C}} is flat on $Z$ 
(or equivalently locally free). We show in \ref{fam_qlfrt} that in this case 
$\phi$ is a morphism to 
\m{V_{E,reg}}.

\sepprop

\begin{subsub} The case of sheaves on ribbons -- Deformations of ribbons -- \rm
We first give a general idea of the way question (iii') can be treated: we 
consider a deformation of $Y$, i.e. a flat morphism \ 
\m{\pi:\kc\to S}, and \m{P\in S}, such that
\begin{enumerate}
\item[--] $S$ is a smooth curve.
\item[--] $\kc_P\simeq Y$.
\item[--] for every $s\in S$, $\kc_s$ is a projective curve.
\end{enumerate}
We suppose that $E$ is a stable quasi locally free sheaf of rigid type. Let \ 
\m{\tau:\M\to S} \ be the {\em relative moduli space of stable sheaves} with 
Hilbert polynomial \m{P_E} (such that for every \m{s\in S}, \m{\M_s} is the 
moduli space of stable sheaves on \m{\kc_s} with Hilbert polynomial \m{P_E},
cf. \cite{ma2}, \cite{si}).
We find that a suitable neighbourhood $U$ of $E$ in \m{\M_P} is deformed in two 
non intersecting open subsets in two components of the fibres \m{\M_s}, 
\m{s\not=P}. The non reduced structure of $U$ comes from the fact that it is 
the ``limit'' of varieties with two connected components.

We have something similar for question (iii). More precisely suppose that \ 
\m{\deg(L)<0}. We will use for $\pi$ a {\em maximal reducible deformation} of 
$Y$, i.e.
\begin{enumerate}
\item[--] $\kc$ is a reduced algebraic variety with two irreducible components 
$\kc_1,\kc_2$.
\item[--] For $i=1,2$, let \ $\pi_i:\kc_i\to S$ \ be the restriction
of $\pi$. Then \ $\pi_i^{-1}(P)=C$ \ and $\pi_i$ is a flat family of smooth
irreducible projective curves.
\item[--] For every $s\in S\backslash\{P\}$, the components
$\kc_{1,s}, \kc_{2,s}$ of $\kc_s$ meet transversally.
\end{enumerate}
For every \m{s\in S\backslash\{P\}}, \m{\kc_{1,s}} and \m{\kc_{2,s}} meet in 
exactly \m{-\deg(L)} points. It is proved in \cite{dr9} that such a deformation 
exists if $L$ can be written as \ \m{L=\ko_C(-P_1-\cdots-P_k)}, for distinct 
points \m{P_1,\ldots,P_k} of $C$. In this case we can construct $\kc$ such that 
if $\kz$ is the closure of the set of intersections points of the components 
\m{\kc_{1,s}}, \m{\kc_{2,s}} of \m{\kc_s}, \m{s\not=P}, we have
\[\kz\cap C \ = \ \{P_1,\ldots,P_k\} \ . \]
In this way we define a one dimensional subspace \ \m{\Delta\subset H^0(L^*)}.

Let $E$ be a quasi locally free sheaf of rigid type on $Y$. Let $\ke$ be a 
coherent sheaf on $\kc$, flat on $S$ and such that \ \m{\ke_{|Y}=E}. Suppose 
that $E$ is locally isomorphic to \m{r\ko_Y\oplus\ko_C} with \m{r>0}. Then we 
from \cite{dr9} there are two possibilities: for every \m{s\in 
S\backslash\{P\}} 
in a neighbourhood of $P$
\begin{enumerate}
\item[1)] $\ke_{s|\kc_{1,s}\backslash\kz}$ is locally free of rank $r$ and  
$\ke_{s|\kc_{2,s}\backslash\kz}$ is locally free of rank $r+1$,
\item[2)] $\ke_{s|\kc_{1,s}\backslash\kz}$ is locally free of rank $r+1$ and  
$\ke_{s|\kc_{2,s}\backslash\kz}$ is locally free of rank $r$,
\end{enumerate}
\end{subsub}
(that is: $\ke$ is of rank $r$ on one of the components of $\kc$ and of rank 
\m{r+1} on the other). So we see that $E$ can be deformed in two distinct ways 
to sheaves on the reduced curves with two components.

\medskip

{\bf Remarks: } {\bf 1 -- } It is proved in \cite{dr9} that given a quasi 
locally free sheaf $E$ on $Y$ (i.e. a sheaf locally isomorphic to a direct sum 
\ \m{a\ko_Y\oplus b\ko_C}), there exists a smooth curve \m{S'}, \m{P'\in S'}, a 
morphism \ \m{f:S'\to S} \ such that \m{f(P')=P}, with non zero tangent map at 
\m{P'}, and a coherent sheaf \m{\ke'} on \m{f^*\kc} flat on \m{S'} and such 
that \m{\ke'_{P'}=E}, i.e. $E$ can be deformed to sheaves on the reduced curves 
with two components.

{\bf 2--} In \cite{dr4} many examples of non-empty moduli spaces of stable 
sheaves containing quasi locally free sheaves of rigid type are given.

\sepprop

\begin{subsub} Koda\"\i ra-Spencer elements -- \rm Let \m{\ke_1}, \m{\ke_2} be 
coherent sheaves on $\kc$, flat on $S$, and such that \ 
\m{\ke_{1,P}=\ke_{2,P}=E}. Then we define 
\[\omega_{\ke_1,\ke_2} \ \in \ \Ext^1_{\ko_Y}(E,E)\]
(cf. \ref{Koda}): on the second neighbourhood \m{Y_2} of $Y$ in $\kc$ we have 
exact sequences
\[0\lra E\lra\ke_{i|Y_2}\lra E\lra 0\]
for \m{i=1,2}, and associated elements \ \m{\sigma_i\in\Ext^1_{\ko_{Y_2}}(E,E)}.
The difference \ \m{\omega_{\ke_1,\ke_2}=\sigma_1-\sigma_2} \ lies in 
\m{\Ext^1_{\ko_Y}(E,E)} (this construction generalises that of the Koda\"\i 
ra-Spencer morphism for sheaves on products \m{X\times Y}).

We then consider two cases

{\bf Case A}
\begin{enumerate} 
\item[--] $\ke^{[1]}_s$ is of rank $r$ on $\kc_{2,s}\backslash\kz$ and of rank
$r+1$ on $\kc_{1,s}\backslash\kz$.
\item[--] $\ke^{[2]}_s$ is of rank $r+1$ on $\kc_{2,s}\backslash\kz$ and of rank
$r$ on $\kc_{1,s}\backslash\kz$.
\end{enumerate}

{\bf Case B}
\begin{enumerate}
\item[--] $\ke^{[1]}_s$ and $\ke^{[2]}_s$ are of rank $r$ on 
$\kc_{2,s}\backslash\kz$ and of rank
$r+1$ on $\kc_{1,s}\backslash\kz$.
\end{enumerate}
Let
\[\phi:\Ext^1_{\ko_Y}(E,E)\lra H^0(\EExt^1_{\ko_Y}(E,E))\simeq H^0(L^*)\]
be the canonical morphism. Recall that its kernel \m{H^1(\EEnd(E))} corresponds 
to good deformations of $E$, or deformations parametrised by a reduced variety 
(cf. 1.1.2). Then we have (theorem \ref{M_theor})

\sepprop

{\bf Theorem : } {\em
{\bf 1 -- } In case A, \m{\phi(\omega_{\ke^{[1]},\ke^{[2]}})} generates 
$\Delta$.

{\bf 2 -- } in case B, we have \ \m{\phi(\omega_{\ke^{[1]},\ke^{[2]}})=0}.}

\sepprop

This means that the ``non-reduced part'' of the deformation of $E$, in \ 
\m{T_{v_0}V_E/T_{v_0}V_{E,reg}}, corresponds to some parameter of the 
deformation of $Y$ to reduced curves with two components.
\sepprop

\begin{subsub}
The case of moduli spaces of stable sheaves -- \rm Let $E$ be a stable quasi 
locally free sheaf of rigid type on $Y$. 
Let \ \m{\tau:\M\to S} \ be the {\em relative moduli space of stable sheaves} 
with Hilbert polynomial \m{P_E} (such that for every \m{s\in S}, \m{\M_s} is 
the moduli space of stable sheaves on \m{\kc_s} with Hilbert polynomial 
\m{P_E}). For \m{s\in S\backslash\{P\}}, let \m{U_{1,s}} (resp. \m{U_{2,s}}) be 
the open subset of \m{\M_s} corresponding to linked sheaves (cf. \ref{C_red}) 
of rank $r$ on \m{\kc_{1,s}} and of rank \m{r+1} on \m{\kc_{2,s}} (resp. of 
rank \m{r+1} on \m{\kc_{1,s}} and of rank $r$ on \m{\kc_{2,s}}).
\end{subsub}

For \m{i=1,2}, let \ \m{{\bf U}_i=\dsp\bigcup_{s\in S\backslash\{P\}}U_{i,s}}.
It is a smooth open subset of $\M$. Let \m{\U_i\subset\M} be the closure of 
\m{{\bf U}_i} and \ \m{\tau_i:\U_i\to S} \ the restriction of $\tau$. Let \ 
\m{\U=\U_1\cup\U_2}, and \ \m{\ov{\tau}:\U\to S} \ the restriction of $\tau$. 
Then from \cite{dr9}, \m{\ov{\tau}^{-1}(P)\subset\M_P} \ contains the points 
corresponding to stable quasi locally free sheaves of rigid type.

\sepprop

The preceding theorem suggests the following natural

{\bf Conjecture: }{\em
{\bf 1 -- } \m{\U_i} is smooth and \m{\tau_i^{-1}(P)} is smooth at $E$, i.e. 
around $E$, we have \ \m{\tau_i^{-1}(P)=(\M_P)_{reg}}.

{\bf 2 -- } \m{\U_1}, \m{\U_2} intersect transversally, and the image of the 
composition
\[T_E(\tau_i^{-1}(P))\lra T_E\M_P=\Ext^1_{\ko_Y}(E,E)\lra 
H^0(\EExt^1_{\ko_Y}(E,E))=H^0(L^*)\] 
is $\Delta$.}

\end{subsub}

\end{sub}

\sepsub

\Ssect{Unstable rank-2 vector bundles on $\P_2$}{P2}

The following results are proved in \cite{st}:

Let $E$ be a rank-2 vector bundle on \m{\P_2} such that \m{c_1(E)=0}. Let 
\m{d(E)} be the largest integer such that \ \m{h^0(E(-d(E)))>0} ($E$ is stable 
if and only \m{d(E)<0}).

Let \m{d,c_2} be integers, such that \m{d\geq 3} and \ \m{d^2+c_2>0}. Let 
\m{B(d)} be the set of isomorphism classes of rank-2 vector bundles $E$ such 
that \m{c_2(E)=c_2} and \m{d(E)=d}.

There is a natural structure of smooth irreducible algebraic variety on 
\m{B(d)}, and \Nligne \m{\dim(B(d))=3(d^2+c_2)-1}. Moreover, the subset $U$ of 
\m{B(d)} consisting of bundles $E$ such that the deformations of $E$ are rank-2 
bundles $F$ such that \m{d(F)=d} is non-empty and open.

We have for \m{E\in B(d)}
\[\dim(\Ext^1_{\ko_{\P_2}}(E,E)) \ = \ 4(d^2+c_2)-1 \ > \ \dim(B(d)) \ . \]
If \m{E\in U}, \m{V_{E,reg}} is smooth of dimension \m{\dim(B(d))}, and \m{V_E} 
is not reduced. 

Let $Z$ be an algebraic variety, \m{z_0\in Z}, and $\kf$ a coherent sheaf on 
\m{Z\times\P_2}, flat on $Z$ and such that \m{\kf_z\in U} for every closed 
point \m{z\in Z}. Let \m{E=\kf_{z_0}}. Then we have a morphism
\[\phi:(Z,z_0)\to (V_E,v_0)\]
such that \ \m{(\phi\times I_{\P_2})^*(\ke_E)\simeq\kf}.

Let \ \m{p_Z:Z\times\P_2\to Z}, \m{p_2:Z\times\P_2\to\P_2} \ be the 
projections, and
\[\kl \ = \ p_{Z*}(\kf\ot p_2^*(\ko_{\P_2}(-d))) \ . \]
For every closed point \m{z\in Z}, \m{\kl_z\ot_{\ko_{Z,z}}\C(z)} \ has 
dimension one. We say that {\em $\kf$ has pure type $d$} if $\kl$ is locally 
free (or equivalently if $\kl$ is flat on $Z$). If $\kl$ is pure of type $d$ 
then $\phi$ is a morphism to \m{V_{E,reg}}, i.e. we have a result which is 
similar to that of 1.1.2.

\end{sub}

\sepsub

\Ssect{Outline of the paper}{outline}

Section 2 contains definitions and properties of {\em primitive multiple 
curves} of any multiplicity, with some particular results in multiplicity 2 
(primitive multiple curves of multiplicity 2 are also called {\em double 
curves} or {\em ribbons}). This section contains also a description of the 
generalisation of the Koda\" ira-Spencer morphism that is used here.

Section 3 is devoted to the study of quasi locally free sheaves of rigid type 
on a primitive multiple curve, and to their deformations. In particular we give 
an answer to question (ii) of 1.1.1. Some results are valid in any 
multiplicity. 

In section 4 we recall some definitions concerning {\em maximal reducible 
deformations} of ribbons, i.e. deformations to reduced curves with two 
components intersecting transversally. We recall also some results about 
deformations of quasi locally free sheaves of rigid type on ribbons to sheaves 
on the reduced curves with two components.

In section 5 we prove the main result of this paper, i.e. the theorem in 1.1.4.
\end{sub}

\sepsec

\section{Preliminaries}

\Ssect{Primitive multiple curves and quasi locally free sheaves}{cmpr}

(cf. \cite{ba_fo}, \cite{ba_ei}, \cite{dr2}, \cite{dr1}, \cite{dr4}, \cite{dr5},
\cite{dr6}, \cite{ei_gr}).

\begin{subsub} Definitions -- \rm Let $C$ be a smooth connected projective 
curve. A {\em multiple curve with support $C$} is a Cohen-Macaulay scheme $Y$ 
such that \m{Y_{red}=C}.

Let $n$ be the smallest integer such that \m{Y=C^{(n-1)}}, \m{C^{(k-1)}}
being the $k$-th infinitesimal neighbourhood of $C$, i.e. \
\m{\ki_{C^{(k-1)}}=\ki_C^{k}} . We have a filtration \ \m{C=C_1\subset
C_2\subset\cdots\subset C_{n}=Y} \ where $C_i$ is the biggest Cohen-Macaulay
subscheme contained in \m{Y\cap C^{(i-1)}}. We call $n$ the {\em multiplicity}
of $Y$.

We say that $Y$ is {\em primitive} if, for every closed point $x$ of $C$,
there exists a smooth surface $S$, containing a neighbourhood of $x$ in $Y$ as
a locally closed subvariety. In this case, \m{L=\ki_C/\ki_{C_2}} is a line
bundle on $C$ and we have \ \m{\ki_{C_j}=\ki_C^j},
\m{\ki_{C_{j}}/\ki_{C_{j+1}}=L^j} \ for \m{1\leq j<n}. We call $L$ the line
bundle on $C$ {\em associated} to $Y$. Let \m{P\in C}. Then there exist
elements $y$, $t$ of \m{m_{S,P}} (the maximal ideal of \m{\ko_{S,P}}) whose
images in \m{m_{S,P}/m_{S,P}^2} form a basis, and such that for \m{1\leq i<n}
we have \ \m{\ki_{C_i,P}=(y^{i})} .

The simplest case is when $Y$ is contained in a smooth surface $S$. Suppose
that $Y$ has multiplicity $n$. Let \m{P\in C} and \m{f\in\ko_{S,P}}  a local
equation of $C$. Then we have \ \m{\ki_{C_i,P}=(f^{i})} \ for \m{1<j\leq n},
in particular \m{I_{Y,P}=(f^n)}, and \ \m{L=\ko_C(-C)} .

For any \m{L\in\Pic(C)}, the {\em trivial primitive curve} of multiplicity $n$,
with induced smooth curve $C$ and associated line bundle $L$ on $C$ is the
$n$-th infinitesimal neighbourhood of $C$, embedded by the zero section in the
dual bundle $L^*$, seen as a surface.

We will write \m{\ko_n=\ko_{C_n}} and we will see \m{\ko_i} as a coherent sheaf
on \m{C_n} with schematic support \m{C_i} if \m{1\leq i<n}.

\end{subsub}

\sepprop

\begin{subsub}\label{inv} Invariants of sheaves and canonical filtrations -- \rm
If $\ke$ is a coherent sheaf on $Y$ one defines its {\em generalised rank
\m{R(\ke)}} and {\em generalised degree \m{\Deg(\ke)}} (cf. \cite{dr4}, 3-): 
take any filtration of $\ke$
\[0=\kf_n\subset\kf_{n-1}\subset\cdots\subset\kf_0=\ke\]
by subsheaves such that \m{\kf_i/\kf_{i+1}} is concentrated on $C$ for 
\m{0\leq i< n}, then
\[R(\ke)=\sigg_{i=0}^{n-1}\text{rk}(\kf_i/\kf_{i+1})  \quad\text{and}\quad
\Deg(\ke)=\sigg_{i=0}^{n-1}\deg(\kf_i/\kf_{i+1}) .\]
Let \m{\ko_Y(1)} be a very ample line bundle on $Y$. Then the Hilbert
polynomial of $\ke$ is
\[P_\ke(m) \ = \ R(\ke)\deg(\ko_C(1))m+\Deg(\ke)+R(\ke)(1-g)\]
(where $g$ is the genus of $C$).

The {\em first canonical filtration} of $\ke$
\[0=\ke_n\subset\ke_{n-1}\subset\cdots\subset\ke_0=\ke\]
is defined as follows: for \m{0\leq i\leq n}, we have \ \m{\ke_i=\ki_C^i\ke}. 
For \m{0\leq i<n}, the sheaf \Nligne \m{G_i(\ke)=\ke_i/\ke_{i+1}} is 
concentrated on $C$. The same definition applies if $\ke$ is a coherent sheaf 
on a non-empty open subset of $Y$. The pair
\[\sigma(\ke) \ = \ \big((\rg(G_0(\ke)),\ldots,\rg(G_{n-1}(\ke)) \ , \
(\deg(G_0(\ke)),\ldots,\deg(G_{n-1}(\ke))\big)\]
is called the {\em complete type} of $\ke$.

Let \m{P\in C} and $M$ a \m{\ko_{Y,P}}-module of finite type. We also define 
the {\em first canonical filtration} of $M$
\[0=M_n\subset M_{n-1}\subset\cdots\subset M_0=M\]
as: \m{M_i=\ki_{C,P}^iM}. The quotients \m{G_i(M)=M_i/M_{i+1}} are 
\m{\ko_{C,P}}-modules. The {\em generalised rank of $M$} is \ 
\m{R(M)=\sigg_{i=0}^{n-1}\rk(G_i(M))}.

Let $X$ be an algebraic variety and $\E$ a coherent sheaf on \m{X\times Y}, 
flat on $X$. We can also define the first canonical filtration of $\E$
\[0=\E_n\subset\E_{n-1}\subset\cdots\subset\E_0=\E\]
by \ \m{\E_i=p_Y^*(\ki_C^i)\E} .

The {\em second canonical filtration} of $\ke$
\[0=\ke^{(0)}\subset \ke^{(1)}\subset\cdots\subset\ke^{(n)}=\ke\]
is defined as follows: for \m{0\leq i\leq n} and \m{P\in C}, \m{\ke^{(i)}_P} is 
the set of \m{u\in\ke_P} such that \ \m{\ki_{C,P}^iu=0}. For \m{1\leq i\leq n}, 
the sheaf \m{G^{(i)}(\ke)=\ke^{(i)}/\ke^{(i-1)}} is concentrated on $C$.

Let \m{P\in C} and $M$ a \m{\ko_{Y,P}}-module of finite type. We also define 
the {\em second canonical filtration} of $M$
\[0=M^{(0)}\subset M^{(1)}\subset\cdots\subset M^{(n)}=M\]
in the obvious way.
\end{subsub}

\sepprop

\begin{subsub} The case of double curves -- \rm If \m{n=2}, let $\ke$ be a 
coherent sheaf on \ \m{Y=C_2}. Then we have canonical exact sequences
\[0\lra\ke_1\lra\ke\lra\ke_{|C}\lra 0 \ , \qquad
0\lra\ke^{(1)}\lra\ke\lra\ke^{(2)}\simeq\ke_1\ot L^*\lra 0 \ . \]
\end{subsub}

\sepprop

\begin{subsub} Quasi locally free sheaves -- \rm Let \m{P\in C} and $M$ a 
\m{\ko_{Y,P}}-module of finite type. We say that $M$ is {\em quasi free} if 
there exist integers \m{m_i\geq 0}, \m{1\leq i\leq n}, such that \ 
\m{M\simeq\oplus_{i=1}^nm_i\ko_{i,P}}. These integers are uniquely determined. 
In this case we say that $M$ is {\em of type} \m{(m_1,\ldots,m_n)}. We 
have \ \m{R(M)=\sigg_{i=1}^ni.m_i}.

Let $\ke$ be a coherent sheaf on a non-empty open subset \m{V\subset Y}. We say 
that $\ke$ est {\em quasi locally free} at a point $P$ of $V$ if there exists a 
neighbourhood \m{U\subset V} of $P$ and integers \m{m_i\geq 0}, \m{1\leq i\leq 
n}, such that for every \m{Q\in U}, \m{\ke_Q} is quasi free of type 
\m{m_1,\ldots,m_n}. The integers \m{m_1,\ldots,m_n} are uniquely determined and 
depend only of $\ke$, and \m{(m_1,\ldots,m_n)} is called the {\em type of} 
$\ke$.

We say that $\ke$ est {\em quasi locally free} if it is quasi locally free at 
every point of $V$.

The following conditions are equivalent: 
\begin{enumerate}
\item[(i)] $\ke_P$ is quasi locally free at $P$.
\item[(ii)] The $\ko_{C,P}$-modules $G_i(\ke_P)$ are free.
\end{enumerate}

The following conditions are equivalent: 
\begin{enumerate}
\item[(i)] $\ke$ is quasi locally free.
\item[(ii)] the sheaves $G_i(\ke)$ are locally free on $C$.
\end{enumerate}

\end{subsub}

\end{sub}

\sepsub

\Ssect{Infinitesimal deformations of coherent sheaves}{Koda1}

\begin{subsub}\label{def_sh} Deformations of sheaves -- \rm Let $X$ be a 
projective algebraic variety and $E$ a coherent sheaf on $X$. A {\em 
deformation} of $E$ is a quadruplet \ \m{\kd=(S,s_0,\ke,\alpha)}, where 
\m{(S,s_0)} is the germ of an analytic variety, $\ke$ is a coherent sheaf on 
\m{S\times X}, flat on $S$, and $\alpha$ an isomorphism \m{\ke_{s_0}\simeq E}.
If there is no risk of confusion, we also say that $\ke$ is an infinitesimal 
deformation of $E$. Let \m{Z_2=\spec(\C|t]/(t^2))}. When \m{S=Z_2} and \m{s_0} 
is the closed point $*$ of \m{Z_2}, we say that $\kd$ is an {\em infinitesimal 
deformation} of $E$. Isomorphisms of deformations of $E$ are defined in an 
obvious way. If \ \m{f:(S',s'_0)\to(S,s_0)} \ is a morphism of germs, the 
deformation \m{f^\#(\kd)} is defined as well. A deformation \ 
\m{\kd=(S,s_0,\ke,\alpha)} \ is called {\em semi-universal} if for every 
deformation \ \m{\kd'=(S',s'_0,\ke',\alpha')} \ of $E$, there exists a morphism 
\ \m{f:(S',s'_0)\to(S,s_0)} \ such that \ \m{f^\#(\kd)\simeq\kd'}, and if the 
tangent map \ \m{T_{s'_0}S'\to T_{s_0}S} \ is uniquely determined. There always 
exists a semi-universal deformation of $E$ (cf. \cite{si_tr}, theorem I).

Let $\ke$ be an infinitesimal deformation of $E$. Let \m{p_X} denote the 
projection \ \m{Z_2\times X\to X}. Then there is a canonical exact sequence
\[0\lra E\lra p_{X*}(\ke)\lra E\lra 0 \ , \]
i.e. an extension of $E$ by itself. In fact, by associating this extension to 
$\ke$ one defines a bijection between the set of isomorphism classes of 
infinitesimal deformations of $E$ and the set of isomorphism classes of 
extensions of $E$ by itself, i.e. \m{\Ext^1_{\ko_X}(E,E)}.
\end{subsub}

\sepprop

\begin{subsub}\label{ko_da} Koda\"\i ra-Spencer morphism -- \rm Let \ 
\m{\kd=(S,s_0,\ke,\alpha)} \ be a deformation of $E$, and \m{X_{s_0}^{(2)}} the 
infinitesimal neighbourhood of order 2 of \ \m{X_{s_0}=\{s_0\}\times X} \ in 
\m{S\times X}. Then we have an exact sequence on \m{X_{s_0}^{(2)}}
\[0\lra T_{s_0}S\ot E\lra\ke/m_s^2\ke=\ke_{|X_{s_0}^{(2)}}\lra E\lra 0 \ . \]
By taking the direct image by \m{p_X} we obtain the exact sequence on $X$
\[0\lra T_{s_0}S\ot E\lra F\lra E\lra 0 \ , \]
hence a linear map
\[\omega_{s_0}:T_{s_0}S\lra\Ext^1_{\ko_X}(E,E) \ , \]
which is called the {\em Koda\"\i ra-Spencer morphism} of $\ke$ at \m{s_0}.

We say that $\ke$ is a {\em complete deformation} if \m{\omega_{s_0}} is 
surjective. If $\kd$ is a semi-universal deformation, \m{\omega_{s_0}} is an 
isomorphism.
\end{subsub}

\end{sub}

\sepsub

\Ssect{Generalisation of the Koda\"\i ra-Spencer morphism}{Koda}

Let $S$ be a smooth curve and \m{s_0\in S} a closed point. Let \ \m{\rho:\kx\to 
S} \ be a flat projective morphism of algebraic varieties. Let 
\m{Y=\rho^{-1}(s_0)}. It is a projective variety. Let $\ke$, \m{\ke'} be 
coherent sheaves on $\kx$, flat on $S$, such that there exists an isomorphism \
\m{\ke_{|Y}\simeq\ke'_{|Y}}. Let \ \m{E=\ke_{|Y}}. Let \m{Y_2} be the second 
infinitesimal neighbourhood of $Y$ in $\kx$. If \m{s_0^{(2)}} is the second 
infinitesimal neighbourhood of \m{s_0} in $S$, we have \ 
\m{Y_2=\rho^{-1}(s_0^{(2)})}. The ideal sheaf \m{\ki_Y} of $Y$ in \m{Y_2} is 
isomorphic to \m{\ko_Y}. We have canonical exact sequences
\xmat{0\ar[r] & E\ot\ki_Y\ar[r]^-i\fleq[d] & \ke_{|Y_2}\ar[r] & E\ar[r] & 0\\ & 
E\fleq[d]\\ 0\ar[r] & E\ot\ki_Y\ar[r]^-{i'} & \ke'_{|Y_2}\ar[r] & E\ar[r] & 0}
(the injectivity of $i$ and \m{i'} follows easily from the flatness of $\ke$, 
\m{\ke'} over $S$). Let \Nligne \m{\sigma,\sigma'\in\Ext^1_{\ko_{Y_2}}(E,E)} \ 
correspond to these extensions.

Let \ \m{0\to E\to\kf\to E\to 0} \ be an exact sequence of coherent sheaves on 
\m{Y_2}. The canonical morphism \ \m{\kf\ot\ki_Y\to\kf} \ induces a morphism \ 
\m{E\ot\ki_Y\to E}, which vanishes if and only if $\kf$ is concentrated on $Y$. 
In this way we get an exact sequence
\xmat{0\ar[r] & \Ext^1_{\ko_Y}(E,E)\ar[r] & \Ext^1_{\ko_{Y_2}}(E,E)\ar[r]^-p & 
\End(E) \ . }
We have \ \m{p(\sigma)=p(\sigma')=I_E}. So we have
\[\omega_{\ke,\ke'} \ = \ \sigma-\sigma' \ \in \ \Ext^1_{\ko_Y}(E,E) \ . 
\]

\sepprop

{\em Connections with the Koda\"\i ra-Spencer morphism -- } Suppose the $\kx$ 
is the trivial family: \m{\kx=Y\times S}. Let
\[\omega_{s_0}:T_{s_0}S\lra\Ext^1_{\ko_Y}(E,E)\]
be the Koda\"\i ra-Spencer morphism of $\ke$. Suppose that \m{\ke'} is the 
trivial family: \m{\ke'=p_Y^*(E)} (where \m{p_Y} is the projection \ \m{Y\times 
S\to S}). The isomorphism \ \m{\ki_Y\simeq\ko_Y} \ is defined by the choice of 
a generator $t$ of the maximal ideal of \m{s_0} in $S$. Let \m{u} the 
associated element of \m{T_{s_0}S}. Then we have \ 
\m{\rho(\sigma-\sigma')=\omega_{s_0}(u)}.

\end{sub}

\sepsec

\section{Quasi locally free sheaves of rigid type}\label{rig_t}

We keep the notations of \ref{cmpr}.

\sepprop

\Ssect{Definitions and basic properties}{QLLRT1}

A quasi locally free sheaf $\ke$ on $Y$ is called {\em of rigid type} if it is 
locally free or locally isomorphic to \m{a\ko_n\oplus\ko_k} for some integers 
\m{a\geq 0}, \m{1\leq k<n}. {\em The set of isomorphism classes of quasi 
locally free sheaves of rigid type of fixed complete type (cf. \ref{inv}) is an 
open family} 
(cf. \cite{dr4}, 6-): let $X$ be an algebraic variety, $\E$ a coherent sheaf 
on \m{X\times Y}, flat on $X$, and \m{x\in X} a closed point. Suppose that 
\m{\E_x} is quasi locally free of rigid type. Then there exists an open subset 
$U$ of $X$ containing $x$ such that, for every \m{x'\in U}, \m{\E_{x'}} is 
quasi locally free and \ \m{\sigma(\E_{x'})=\sigma(\E_x)}.

More generally, let $\ke$ be a quasi locally free sheaf on $Y$, locally 
isomorphic to \m{a\ko_n\oplus b\ko_k}, with \m{a,b>0}, \m{1\leq k<n}. By 
\cite{dr4}, prop. 5.1, there exists a vector bundle $\E$ on $Y$ and a 
surjective morphism
\[\phi:\E\lra\ke\]
inducing an isomorphism \ \m{\E_{|C}\simeq\ke_{|C}}. Let \ \m{\kl=\ker(\phi)}. 
By \cite{dr4}, lemme 5.2, $\kl$ is a vector bundle of rank $b$ on \m{C_{n-k}}.
Let \m{P\in C}, and \m{z\in\ko_{Y,P}} an equation of $C$.
At $P$ the exact sequence \ \m{0\to\kl\to\E\to\ke\to 0} \ is isomorphic to the 
trivial one 
\[0\lra(z^k)\ot\C^b\simeq\ko_{n-k,P}\ot\C^b\lra\ko_{n,P}\ot(\C^a\oplus\C^b)\lra
(\ko_{n,P}\ot\C^a)\oplus(\ko_{k,P}\ot\C^b)\lra 0 . \ \]

\sepprop

\begin{subsub}\label{lem1}{\bf Lemma: } There is a canonical isomorphism
\[\kl_{|C_k} \ \simeq \ (\ke^{(k)}/\ke_{n-k})\ot\ki_C^k \ . \]
\end{subsub}

\begin{proof}
Let \m{P\in C}, \m{z\in\ko_{Y,P}} an equation of $C$ and \m{u\in\ke^{(k)}_P}.
Let \m{v\in\E_P} be such that \m{\phi_P(v)=u}. Then we have
\[\phi_P(z^kv) \ = \ z^ku \ = \ 0 \ , \]
hence \m{z^kv\in\kl_P}. If \m{v'\in\E_P} is such that \m{\phi_P(v')=u}, we have 
\m{w=v'-v\in\kl_P}, hence the image of \m{z^kv'=z^kv+z^kw} in 
\m{\kl_P/z^k\kl_P} is the same as that of \m{z^kv}. By associating \m{z^kv} to 
\m{u\ot z^k} we define a morphism \ 
\m{\ov{\theta}:\ke^{(k)}\ot\ki_C^k\to\kl_{|C_k}}. If \m{u\in\ke_{n-k,P}}, let 
\m{u'\in\ke_P} be such that \m{u=z^{n-k}u'}. Let \m{v'\in\E_P} be such that 
\m{\phi_P(v')=u'}. Then we can take \m{v=z^{n-k}v'}, and then \m{z^kv=0}, hence 
\m{\ov{\theta}_P(u\ot z^k)=0}. It follows that $\ov{\phi}$ induces a morphism
\[\theta:(\ke^{(k)}/\ke_{n-k})\ot\ki_C^k\lra\kl_{|C_k} \ . \]

In the above description of the exact sequence \ \m{0\to\kl\to\E\to\ke\to 0}, 
\Nligne we have \ \m{\ke^{(k)}_P/\ke_{n-k,P}=\ko_{k,P}\ot\C^b}, and \ 
\m{\theta_P=(\ko_{k,P}\ot\C^b)\ot(z^k)\to(z^k)\ot\C^b} is the identity 
morphism. Hence $\phi$ is an isomorphism.
\end{proof}

\sepprop

The sheaf \m{\kl_{|C_k}} is a vector bundle of rank $b$, on \m{C_k} if 
\m{2k\leq n}, and on \m{C_{n-k}} if \m{2k>n}.

\sepprop

\begin{subsub}\label{coro1}{\bf Corollary: } There is a canonical isomorphism
\[\EExt^1_{\ko_Y}(\ke,\ke) \ \simeq \
\HHom((\ke^{(k)}/\ke_{n-k})\ot\ki_C^k,\ke^{(k)}/\ke_{n-k}) \ . \]
\end{subsub}

\begin{proof} From the exact sequence \ \m{0\to\kl\to\E\to\ke\to 0}, we deduce
the exact sequence
\[\HHom(\E,\ke)\lra\HHom(\kl,\ke)\lra\EExt^1_{\ko_Y}(\ke,\ke)\lra 0 \ , \]
and the result follows easily using local isomorphisms of $\ke$ with 
\m{a\ko_n\oplus b\ko_k} and lemma \ref{lem1}.
\end{proof}

\sepprop

If $\ke$ is of rigid type (i.e. if \m{b=1}), then \m{\ke^{(k)}/\ke_{n-k}} is a 
line bundle on \m{C_{\inf(k,n-k)}}, and it follows that
\[\EExt^1_{\ko_Y}(\ke,\ke) \ \simeq \ 
\HHom(\ki_{C|C_{\inf(k,n-k)}}^k,\ko_{\inf(k,n-k)}) \ . \]
It follows that we have an exact sequence
\xmat{0\ar[r] & H^1(\EEnd(\ke))\ar[r] & \Ext^1_{\ko_Y}(\ke,\ke)\ar[r] & 
\Hom(\ki_{C|C_{\inf(k,n-k)}}^k,\ko_{\inf(k,n-k)})\fleq[d]\ar[r] & 0\\
& & & H^0(\EExt^1_{\ko_Y}(\ke,\ke))}

\sepprop

\begin{subsub}\label{coro1b} The case of double curves -- \rm
If \m{n=2} we have \m{k=1} and from corollary \ref{coro1}
\[\EExt^1_{\ko_Y}(\ke,\ke) \ \simeq \ L^* \ . \]
\end{subsub}

\sepprop

\begin{subsub}{\bf Remark: }\rm
Let Let $X$ be an algebraic variety, $\E$ a coherent sheaf on \m{X\times Y}, 
flat on $X$, such that for every closed point \m{x\in X}, \m{\E_x} is quasi 
locally free of rigid type. Let \m{x\in X}, and
\[\omega_x(\E):TX_x\lra\Ext^1_{\ko_Y}(\E_x,\E_x)\]
the Koda\"\i ra-Spencer morphism of $\E$. Let \m{X_\text{red}} be the reduced 
subscheme associated to $X$. Then the image of Koda\"\i ra-Spencer morphism of 
\m{\E_{X_\text{red}\times Y}}
\[\omega_x(\E_{X_\text{red}\times 
Y}):TX_{\text{red},x}\lra\Ext^1_{\ko_Y}(\E_x,\E_x)\]
is contained in \m{H^1(\EEnd(\E_x))}. Suppose that $\E$ is a complete 
deformation of \m{\E_x} (i.e. \m{\omega_x(\E)} is surjective), and that 
\m{\E_x} is simple. Then \ 
\m{\imm(\omega_x(\E_{|X_\text{red}\times Y}))=H^1(\EEnd(\E_x))} (cf. \cite{dr4}, 
th. 6.10, cor. 6.11).
\end{subsub}

\end{sub}

\sepsub

\Ssect{Families of quasi locally free sheaves of rigid type}{fam_qlfrt}

Let \m{m_1,\ldots,m_n} be non negative integers. Let $X$ be a 
connected algebraic variety, \m{U\subset X\times Y} an open subset such that \ 
\m{p_X(U)=X} (where \m{p_X} is the projection \m{X\times Y\to X}) and $\E$ a 
coherent sheaf on $U$, flat on $X$, such that for every closed point 
\m{x\in X}, \m{\E_x} is quasi locally 
free of type \m{(m_1,\ldots,m_n)}. We say that $\E$ is a {\em good family} if 
for \m{0\leq i<n} the sheaf \m{\E_i/\E_{i+1}} on \m{(X\times C)\cap U} is flat 
on $X$ (where \m{0=\E_n\subset\E_{n-1}\subset\cdots\subset\E_0=\E} is the first 
canonical filtration of $\E$). It $\E$ is a good family then by \cite{sga1}, 
exp. IV, prop. 1.1, for \m{0\leq i<n}, \m{\E_i} is a flat family of sheaves on 
\m{C_{n-i}}, and by \cite{si}, lemma 1.27, \m{\E_i/\E_{i+1}} is a vector bundle 
on \m{X\times C}.

\sepprop

\begin{subsub}\label{theo1}{\bf Theorem: } {\bf 1 -- } The sheaf $\E$ is a good 
family if and only if it is locally isomorphic to \ 
\m{\bigoplus_{i=1}^n\ko_{X\times C_i}\ot\C^{m_i}}. 

If $\E$ is a good family on \m{X\times Y}, then for every \m{x\in X} the image 
of Koda\"\i ra-Spencer morphism of $\E$
\[\omega_x(\E):TX_x\lra\Ext^1_{\ko_Y}(\E_x,\E_x)\]
is contained in \m{H^1(\EEnd(\E_x))}.
\end{subsub}

\begin{proof} Suppose that $\E$ is locally isomorphic to \ 
\m{\bigoplus_{i=1}^n\ko_{X\times C_i}\ot\C^{m_i}}. Then it is obvious that the 
sheaves \m{\E_i/\E_{i+1}} are vector bundles on \m{(X\times C)\cap U}, hence 
they are flat on $X$ and $\E$ is a good family. On the other hand, the local 
structure of $\E$ does not vary when $x$ varies, hence for every \m{P\in C}, 
the image of \m{\imm(\omega_x(\E))} in 
\m{\Ext^1_{\ko_{Y,P}}(\E_{x,P},\E_{x,P})} must be 0, so \ 
\m{\imm(\omega_x(\E))\subset H^1(\EEnd(\E_x))}.

Conversely, suppose that $\E$ is a good family. The proof that $\E$ is locally 
isomorphic to \ \m{\bigoplus_{i=1}^n\ko_{X\times C_i}\ot\C^{m_i}} is similar to 
that of theor. 6.5 of \cite{dr4}. We make an induction on $n$. The result for 
\m{n=1} follows from \cite{si}, lemma 1.27. Suppose that it is true for 
\m{n-1\geq 1}. We make an induction on \m{m_n}.

Suppose that \m{m_n=0}. Let $k$ be the smallest integer such that \m{m_q=0} for 
\m{k+1\leq q\leq n}. Then we have \m{k<n}, and for every \m{x\in X}, \m{\E_x} 
is concentrated on \m{C_k}. Then $\E$ is concentrated on \m{(X\times C_k)\cap 
U}: this follows easily by induction on $k$ from the exact sequence \ 
\m{0\to\E_1\to\E\to\E_{(X\times C)\cap U}\to 0}, using the fact that \m{\E_1} 
is flat on $X$. By the induction hypothesis (on \m{(X\times C_k)\cap U}), $\E$ 
is locally isomorphic to \ \m{\bigoplus_{i=1}^n\ko_{X\times C_i}\ot\C^{m_i}}.

Suppose that the result is true for \m{m_n-1\geq 0}. Let \m{P\in C}, \m{x\in X} 
such that \m{Q=(x,P)\in U}. For every open subset \m{V\subset X\times Y} and 
\m{x'\in X}, let \m{V_{x'}=V\cap(\{x\}\times Y)}. Let \m{Z\subset U} be an open 
affine subset containing $Q$ such that there is an isomorphism
\[\E_{x|Z_x} \ \simeq \ \bigoplus_{i=1}^n\ko_{i|Z_x}\otimes\C^{m_i} \ , \]
and that \ \m{\Delta=p_Y^*(\ki_C)} (which is a line bundle on \m{(X\times 
C_{n-1})\cap U}) is trivial on \m{(X\times C_{n-1})\cap Z}. Let 
\m{\zeta\in H^0(\Delta)} a section inducing an isomorphism 
\m{\Delta\simeq\ko_{(X\times C_{n-1})\cap Z}}. Let \ \m{\sigma\in 
H^0(\E_{x|Z_x})} defined by some non zero element of \m{\C^{m_n}}, and 
\m{\ov{\sigma}\in H^0(\E_{|Z})} extending $\sigma$. Then \ 
\m{\zeta^{n-1}\ov{\sigma}\in H^0(Z,\E_{n-1})}, \m{\E_{n-1}} is a vector 
bundle on \m{(X\times C)\cap Z}, and \ \m{s=\zeta^{n-1}\ov{\sigma}_{|Z_x}} \ 
does not vanish on \m{Z_x}. Let \m{T\subset Z} be the open subset where $s$ 
does not vanish. Let \m{x'\in X} such that \m{T_{x'}\not=\emptyset}, and 
\m{W\subset T_{x'}} an open subset such that \ 
\m{\E_{|W}\simeq\bigoplus_{i=1}^n\ko_{i|W}\times\C^{m_i}}. Then in \ 
\m{\ov{\sigma}_{|W}:\ko_{n|W}\to\ko_{n|W}\to\C^{m_n}} \ does not vanish at any 
point. It follows that
\[\coker(\ov{\sigma}_{|W}) \ \simeq \ \left(\bigoplus_{i=1}^{n-1}\ko_{i|W}\ot
\C^{m_i}\right)\oplus\left(\ko_{n|W}\ot\C^{m_n-1}\right) \ . \]
From \cite{sga1}, exp. IV, cor. 5.7, \m{\kf=\coker(\ov{\sigma}_{|T})} is flat on 
$X$. It is a family of quasi locally free sheaves of type 
\m{(m_1,\ldots,m_{n-1},m_n-1)}, and it is easy to verify that it is a good 
family. From the induction hypothesis we can assume, by replacing $T$ with a 
smaller affine neighbourhood of $Q$, that
\[\kf \ \simeq \ \left(\bigoplus_{i=1}^{n-1}\ko_{(X\times C_i)\cap T}\ot
\C^{m_i}\right)\oplus\left(\ko_T\ot\C^{m_n-1}\right) \ . \]
Hence we have an exact sequence
\[0\lra\ko_T\lra\E_{|T}\lra\left(\bigoplus_{i=1}^{n-1}\ko_{(X\times C_i)\cap 
T}\ot\C^{m_i}\right)\oplus\left(\ko_T\ot\C^{m_n-1}\right)\lra 0 \ . \]
Now we have \ \m{\Ext^1_{\ko_T}(\ko_{(X\times C_i)\cap T},\ko_T)=\nsp} for 
\m{1\leq i\leq n} : it suffices to prove that \Nligne
\m{\EExt^1_{\ko_T}(\ko_{(X\times C_i)\cap T},\ko_T)=0}.
This follows easily from the resolution
\xmat{\cdots\ar[r] & \ko_T\ar[r]^{\times \zeta^i} & \ko_T\ar[r]^{\times
\zeta^{n-i}} & \ko_T\ar[r]^{\times \zeta^i} & \ko_T\ar[r] & \ko_{(X\times 
C_i)\cap T} \ . }
Hence
\[\E_{|T} \ \simeq \bigoplus_{i=1}^n\ko_{(X\times C_i)\cap T}\ot\C^{m_i} \ , \]
and the result is proved for \m{m_n}.
\end{proof}

\sepprop

Let $\E$ be a good family of quasi locally free sheaves of rigid type 
parametrised by $X$, and \m{x\in X} a closed point. Suppose that \m{\E_x} is 
simple. Let \m{(S,s_0,\ke,\alpha)} be a semi-universal deformation of \m{\E_x} 
(cf. \ref{def_sh}). Let \m{f:S(x)\to S} be the morphism induced by $\E$ (where 
\m{S(x)} is the germ defined by $\E$ around $x$). Then \m{TS_{s_0}} is 
canonically isomorphic to \m{\Ext^1_{\ko_Y}(\E_x,\E_x)}, and by \cite{dr4}, th. 
6.10 and cor. 6.11, we have
\[TS_{\text{reg},s_0} \ = \ H^1(\EEnd(\E_x)) \ . \]
It follows easily from theorem \ref{theo1} that the image of \ 
\m{Tf_{x}:TX_x\to TS_{s_0}} \ is contained in \m{TS_{\text{reg},s_0}} and that 
the image of $f$ is contained in \m{S_\text{reg}}. Hence if $\boldsymbol{M}$ is 
the moduli space of stable sheaves corresponding to \m{\E_x} and $X$ is 
connected, the image of the canonical morphism \ \m{f_\E:X\to\boldsymbol{M}} 
associated to $\E$ is contained in \m{\boldsymbol{M}_\text{reg}}.

\end{sub}

\sepsec

\section{Coherent sheaves on reducible deformations of primitive double curves}

\Ssect{Maximal reducible deformations}{prel2}

(cf. \cite{dr7}, \cite{dr8}, \cite{dr9})

Let $C$ be a projective irreducible smooth curve and \m{Y=C_2} a primitive 
double curve, with underlying smooth curve $C$, and associated line bundle $L$ 
on $C$. Let $S$ be a smooth curve, \m{P\in S} and \ \m{\pi:\kc\to S} \ a {\em 
maximal reducible deformation of \m{Y}} (cf. \cite{dr7}). This means that
\begin{enumerate}
\item[(i)] $\kc$ is a reduced algebraic variety with two irreducible components 
$\kc_1,\kc_2$.
\item[(ii)] We have \ $\pi^{-1}(P)=Y$. So we can view $C$ as a curve in
$\kc$.
\item[(iii)] For $i=1,2$, let \ $\pi_i:\kc_i\to S$ \ be the restriction
of $\pi$. Then \ $\pi_i^{-1}(P)=C$ \ and $\pi_i$ is a flat family of smooth
irreducible projective curves.
\item[(iv)] For every $z\in S\backslash\{P\}$, the components
$\kc_{1,z},\kc_{2,z}$ of $\kc_z$ meet transversally.
\end{enumerate}

For every \m{z\in S\backslash\{P\}}, \m{\kc_{1,z}} and \m{\kc_{2,z}} meet in 
exactly \m{-\deg(L)} points. If \ \m{\deg(L)=0}, then $\pi$ (or $\kc$) is 
called a {\em fragmented deformation}.

Let \m{\kz\subset\kc} be the closure in $\kc$ of the locus of the intersection
points of the components of \m{\pi^{-1}(z)}, \m{z\not=P}. Since $S$ is a curve,
$\kz$ is a curve of \m{\kc_1} and \m{\kc_2}. It intersects $C$ in a finite
number of points. If \m{x\in C}, let \m{r_x} be the number of
branches of $\kz$ at $x$ and \m{s_x} the sum of the multiplicities of the
intersections of these branches with $C$. If \m{x\in\kz}, then the branches of 
$\kz$ at $x$ intersect transversally with $C$, and we have \m{r_x=s_x}. We have
\[L \ \simeq \ \ko_C(-\sigg_{x\in \kz\cap C}r_xx) \ \simeq \ \ki_{\kz\cap C,C}
\ . \]
For every \m{\in C }, there exists an unique integer \m{p>0} such that 
\m{\ki_{C,x}/\span{(\pi_1,\pi_2)}} is generated by the image of 
\m{(\pi_1^p\lambda_1,0)}, for some \m{\lambda_1\in\ko_{\kc_1,x}} not divisible 
by \m{\pi_1}. Moreover \m{(\pi_1^p\lambda_1,0)} is a generator of the ideal
\m{\ki_{\kc_1,\kc,x}} of \m{\kc_2} in $\kc$, and \m{\lambda_1} is a generator 
of the ideal of $\kz$ in \m{\kc_1} at $x$. The integer $p$ does not depend on 
$x$. Of course we have a symmetric result: \m{\ki_{C,x}/\span{(\pi_1,\pi_2)}} 
is generated by the image of  \m{(0,\pi_2^p\lambda_2)}, for some 
\m{\lambda_2\in\ko_{\kc_2,x}} not divisible by \m{\pi_2}. Moreover 
\m{(0,\pi_2^p\lambda_2)} is a generator of the ideal \m{\ki_{\kc_2,\kc,x}} of 
\m{\kc_1} in $\kc$, and \m{\lambda_2} is a generator of the ideal of $\kz$ in 
\m{\kc_2} at $x$. We can even assume that \ 
\m{(\lambda_1,\lambda_2)\in\ko_{\kc,x}}.

Let
\[\kz_0 \ = \ \kc_1\cap\kc_2 \ \subset \kc \ . \]
We have then \ \m{(\kz_0)_{red}=\kz\cup C} (and \m{\kz_0=\kz\cup C} if 
\m{p=1}). The ideal sheaf \ 
\m{\L_1=\ki_{Z_0,\kc_1}} (resp. \m{\L_2=\ki_{Z_0,\kc_2}}) of \m{\kz_0} in 
\m{\kc_1} (resp. \m{\kc_2}) at $x$ is generated by \m{\lambda_1\pi_1^p} (resp. 
\m{\lambda_2\pi_2^p}). Hence \m{\L_1} (resp. \m{\L_2}) is a line bundle on 
\m{\kc_1} (resp. \m{\kc_2}). The ideal sheaf  \ \m{\ki_{Z,\kc_1}} (resp. 
\m{\ki_{Z,\kc_2}}) of $\kz$ in \m{\kc_1} (resp. \m{\kc_2}) is canonically 
isomorphic to \m{\L_1} (resp. \m{\L_2}). The $p$-th infinitesimal 
neighbourhoods of $C$ in \m{\kc_1}, \m{\kc_2} (generated respectively by 
\m{\pi_1^p} and \m{\pi_2^p}) are canonically isomorphic, we will denote them by 
\m{C^{(p)}}. We have also a canonical isomorphism \ 
\m{\L_{1|C^{(p)}}\simeq\L_{2|C^{(p)}}}, and \ \m{\L_{1|C}\simeq\L_{2|C}\simeq 
L}. It is also possible, by replacing $S$ with a smaller neighbourhood of $P$, 
to assume that \ \m{\L_{1|\kz_0}\simeq\L_{2|\kz_0}}. Let \ 
\m{\L=\L_{1|\kz_0}=\L_{2|\kz_0}}.

We have \ \m{\ki_{\kc_1,\kc}=\L_2} \ and \ \m{\ki_{\kc_2,\kc}=\L_1}.

In this paper we will always assume that \m{p=1}. 

There exists a maximal reducible deformation of \m{Y} either if \m{\deg(L)=0}, 
or if \m{\deg(L)<0} and there exists \m{-\deg(L)} distinct points 
\m{P_1,\ldots,P_d} of $C$ (with \ \m{d=-\deg(L)}) such that \ 
\m{L=\ko_C(-P_1-\cdots-P_d)}. And in the second case we can even assume that 
\Nligne \m{\kz\cap C=\{P_1,\ldots,P_d\}}.

\end{sub}

\sepsub

\Ssect{Coherent sheaves on reduced reducible curves}{C_red}

(cf. \cite{dr9}, 4-)

Let $D$ be a projective curve with two components \m{D_1}, \m{D_2} intersecting 
transversally, and \ \m{Z=D_1\cap D_2}. Let $\ke$ be a coherent sheaf on $D$. 
Then the following conditions are equivalent:
\begin{enumerate}
\item[(i)] $\ke$ is pure of dimension 1.
\item[(ii)] $\ke$ is of depth 1.
\item[(iii)] $\ke$ is locally free at every point of $X$ belonging to only one 
component, and if $x\in Z$, then 
there exist integers $a$, $a_1$, $a_2\geq 0$ and an isomorphism
\[\ke_x \ \simeq \ a\ko_{X,x}\oplus a_1\ko_{D_1,x}\oplus a_2\ko_{D_2,x} \ . \]
\item[(iv)] $\ke$ is torsion free, i.e. for every $x\in X$, every element of 
$\ko_{X,x}$ which is not a zero divisor in $\ko_{X,x}$ is not a zero divisor in 
$\ke_x$.
\item[(v)] $\ke$ is reflexive.
\end{enumerate}

Let \m{E_i=\ke_{|D_i}/T_i}, where \m{T_i} is the torsion subsheaf. It is a 
vector bundle on \m{D_i}.
Let \m{x\in Z}. Then there exists a finite dimensional vector space $W$, 
surjective maps \ \m{f_i:E_{i,x}\to W}, such that the \m{\ko_{D,x}}-module 
\m{\ke_x} is isomorphic to \ \m{\{(\phi_1,\phi_2)\in E_1(x)\times E_2(x);
f_1(\phi_1(x))=f_2(\phi_2(x))\}} (where \m{E_i(x)} is the fibre at $x$ of the 
sheaf \m{E_i}, and \m{E_{i,x}} the fibre of the corresponding vector bundle). 
We have then
\[\ke_x \ \simeq \ (W\times\ko_{D,x})\oplus(\ker(f_1)\ot\ko_{D_1,x})
\oplus(\ker(f_2)\ot\ko_{D_2,x}) \ . \]
We say that the sheaf is {\em linked} at $x$ if $W$ has the maximal possible 
dimension, i.e. \ \m{\dim(W)=\inf(rk(E_1),\rk(E_2))} (i.e. if in (iii) 
\m{a_1=0} or \m{a_2=0}). We say that $\ke$ is {\em linked} if it is linked at 
every point of $Z$.
\end{sub}

\sepsub

\Ssect{Regular sheaves}{def_rig}

A coherent sheaf $\ke$ 
on $\kc$ is called {\em regular} if it is locally free on 
\m{\kc\backslash\kz_0}, and if for every \m{x\in\kz_0} there exists a 
neighbourhood of $x$ in $\kc$, a vector bundle $\E$ on $U$, \m{i\in\{1,2\}}, 
and a vector bundle $F$ on \m{U\cap\kc_i}, such that \ 
\m{\ke_{|U}\simeq\E\oplus F}.

Let $\ke$ be a coherent sheaf on 
$\kc$. Then by \cite{dr9}, prop. 6.4.3, the following assertions are equivalent:
\begin{enumerate}
\item[(i)] $\ke$ is regular (with $i=1$).
\item[(ii)] There exists an exact sequence \ \m{0\to E_2\to\ke\to E_1\to 0}, 
where for \m{j=1,2}, \m{E_j} is a vector bundle on \m{\kc_j}, such that the 
associated morphism \ \m{E_{1|\kz_0}\to\L^*\ot E_{2|\kz_0}} \ is surjective on
a neighbourhood of $C$.
\item[(iii)] There exists an exact sequence \ \m{0\to E_1\to\ke\to E_2\to 0}, 
where for \m{j=1,2}, \m{E_j} is a vector bundle on \m{\kc_j}, such that the 
associated morphism \ \m{E_{2|\kz_0}\to\L^*\ot E_{1|\kz_0}} \ is injective (as 
a morphism of vector bundles) on a neighbourhood of $C$.
\end{enumerate}
If we restrict the exact sequence of (ii) to $Y$ we get the canonical one
\[0\lra(\ke_{|Y})_1\lra\ke_{|Y}\lra(\ke_{|Y})_{|C}\lra 0 \ , \]
(cf. \ref{cmpr}) and if we restrict the exact sequence of (iii) to $Y$ we get
\[0\lra(\ke_{|Y})^{(1)}\lra\ke_{|Y}\lra(\ke_{|Y})^{(2)}=(\ke_{|Y})_1\ot L^*\lra 
0 \ . \]
In particular \m{\ke_{|Y}} is quasi locally free, and for \m{s\in\backslash\{
P\}} in a neighbourhood of $P$, \m{\ke_s} is a linked torsion free sheaf.

We have a similar result by taking \m{i=2}.

For example, let $\ke$ be a coherent sheaf
on $\kc$, flat on $S$. Suppose that for every \m{s\in S}, \m{\ke_s} is torsion
free, and that \m{\ke_{|Y}} is quasi locally free of rigid type (cf.
\ref{rig_t}). Then $\ke$ is regular (\cite{dr9}, prop. 6.4.5).

\end{sub}

\sepsub

\newpage

\section{Koda\" ira-Spencer elements}

We keep the notations of \ref{prel2}.

\sepprop

\Ssect{Self-extensions of $\ko_{C,x}$ on $Y$}{self}

We will need in \ref{main} a description of the extensions
\begin{equation}\label{equ6}
0\lra\ko_{C,x}\lra\ke\lra\ko_{X,x}\lra 0
\end{equation}
on $Y$.

Let \m{x\in C}. Let \m{z\in\ko_{Y,x}} be an equation of $C$ and 
\m{t\in\ko_{Y,x}} over a generator of the maximal ideal of \m{\ko_{C,x}}. The 
extensions $(\ref{equ6})$ are parametrised by \ 
\m{\Ext^1_{\ko_{Y,x}}(\ko_{C,x},\ko_{C,x})}, which is isomorphic to 
\m{\ko_{C,x}}. This can be seen easily by using the free resolution of 
\m{\ko_{C,x}} on $Y$:
\xmat{\cdots\ar[r] & \ko_{Y,x}\ar[r]^-{\times z} & \ko_{Y,x}\ar[r]^-{\times z} &
\ko_{C,x}\ar[r] & 0}
For every positive integer $n$, let
\[\ki_{Y,n} \ = \ (z,t^n) \ \subset \ \ko_{Y,x} \ , \qquad
\ki_{C,n} \ = \ (t^n) \ \subset \ \ko_{C,x}\]
(the ideals of \m{nx}). Then we have an obvious extension
\xmat{0\ar[r] & (z)\fleq[d]\ar[r] & \ki_{Y,n}\ar[r] & \ki_{C,n}\fleq[d]\ar[r] & 
0\\ & \ko_{C,x} & & \ko_{C,x}}
and it is easy to see that it is associated to \ 
\m{t^n\in\Ext^1_{\ko_{Y,x}}(\ko_{C,x},\ko_{C,x})}.

\end{sub}

\sepsub

\Ssect{Proof of the main result}{main}

Let \m{\ke^{[1]}}, \m{\ke^{[2]}} be coherent sheaves on $\kc$, flat on $S$. 
Suppose that \m{\ke^{[1]}_{|Y}}, \m{\ke^{[2]}_{|Y}} are isomorphic. Let
\[E \ = \ \ke^{[1]}_{|Y} \ = \ \ke^{[2]}_{|Y} \ . \]
Suppose that $E$ is quasi locally free of rigid type, and that for 
every \m{s\in S}, \m{\ke^{[1]}_s} and \m{\ke^{[2]}_s} are torsion free. Then 
\m{\ke^{[1]}} and \m{\ke^{[2]}} are regular (cf. \ref{def_rig}). Suppose that 
for every \m{s\in S\backslash\{P\}}, \m{\ke^{[1]}_s} and \m{\ke^{[2]}_s} are 
linked (this is always true on a neighbourhood of $P$). It follows that there 
exists an integer $r$ such that for \m{i=1,2}, for every \m{s\in 
S\backslash\{P\}}, \m{\ke^{[i]}_s} is of rank $r$ on \m{\kc_{1,s}\backslash\kz} 
and \m{r+1} on \m{\kc_{2,s}\backslash\kz}, or 
of rank $r$ on \m{\kc_{2,s}\backslash\kz} and \m{r+1} on 
\m{\kc_{1,s}\backslash\kz}. We suppose that \m{r>0}. We will consider two cases:

{\bf Case A}
\begin{enumerate} 
\item[--] $\ke^{[1]}_s$ is of rank $r$ on $\kc_{2,s}\backslash\kz$ and of rank
$r+1$ on $\kc_{1,s}\backslash\kz$.
\item[--] $\ke^{[2]}_s$ is of rank $r+1$ on $\kc_{2,s}\backslash\kz$ and of rank
$r$ on $\kc_{1,s}\backslash\kz$.
\end{enumerate}

{\bf Case B}
\begin{enumerate}
\item[--] $\ke^{[1]}_s$ and $\ke^{[2]}_s$ are of rank $r$ on 
$\kc_{2,s}\backslash\kz$ and of rank
$r+1$ on $\kc_{1,s}\backslash\kz$.
\end{enumerate}

\sepprop

We want to study \ \m{\omega_{\ke^{[1]},\ke^{[2]}}\in\Ext^1_{\ko_Y}(E,E)} (cf. 
\ref{Koda}). 

Recall that \ \m{\EExt^1_{\ko_Y}(E,E)\simeq 
H^0(L^*)} (cf. \ref{coro1}, \ref{coro1b}). From \ref{prel2}, $\kc$ induces a 
one dimensional subspace \ \m{\Delta\subset H^0(L^*)}.

Let \m{Y_2} be the second infinitesimal neighbourhood of $Y$ in $\kc$. Let 
\m{t\in\ko_{S,P}} be a generator of the maximal ideal. We will also denote by 
$\pi$ (resp. \m{\pi_i}, \m{i=1,2}) the regular function \m{t\circ\pi} (resp. 
\m{t\circ\pi_i}) defined on a neighbourhood of $C$. Then \m{Y_2} is defined in 
a neighbourhood of $Y$ by the equation \ \m{\pi^2=0}. Let $\ki$ be the 
ideal sheaf of $Y$ in \m{Y_2}. We have \m{\ki\simeq\ko_Y}. For \m{i=1,2} we 
have a canonical exact sequence
\begin{equation}\label{equ1}0\lra E\ot\ki\simeq E\lra\ke^{[i]}_{|Y_2}\lra E\lra 
0 \ , \end{equation}
associated to \ \m{\sigma_i\in\Ext^1_{\ko_{Y_2}}(E\ot\ki,E)}.

Given an extension \ \m{0\to E\ot\ki\to\kf\to E\to 0} \ on \m{Y_2}, the 
canonical morphism \ \m{\kf\ot\ki\to\kf} \ induces an endomorphism of $E$. In 
this way we get a canonical morphism \Nligne
\m{\Ext^1_{\ko_{Y_2}}(E,E\ot\ki)\to\End(E)}, whose kernel corresponds to 
extensions such that $\kf$ is concentrated on $Y$. Hence we have an exact 
sequence
\xmat{0\ar[r] & \Ext^1_{\ko_Y}(E\ot\ki,E)\ar[r] & 
\Ext^1_{\ko_{Y_2}}(E\ot\ki,E)\ar[r]^-\theta & \End(E) \ .}
The image of \m{\sigma_i}, \m{i=1,2}, is \m{I_E}. Hence, by using the action of 
\m{\Aut(E)}, we see that $\theta$ is surjective, and that we have an 
exact sequence
\xmat{0\ar[r] & \Ext^1_{\ko_Y}(E\ot\ki,E)\ar[r] & 
\Ext^1_{\ko_{Y_2}}(E\ot\ki,E)\ar[r]^-\theta & \End(E)\ar[r] & 0 \ .}
Recall that
\[\omega_{\ke^{[1]},\ke^{[2]}} \ = \ \sigma_1-\sigma_2 \ . \]
Let
\[\phi:\Ext^1_{\ko_Y}(E,E)\lra H^0(\EExt^1_{\ko_Y}(E,E))\]
be the canonical morphism. 
 
\sepprop

\begin{subsub}\label{M_theor}{\bf Theorem : } 
{\bf 1 -- } In case A, \m{\phi(\omega_{\ke^{[1]},\ke^{[2]}})} generates 
$\Delta$.

{\bf 2 -- } in case B, we have \ \m{\phi(\omega_{\ke^{[1]},\ke^{[2]}})=0}.
\end{subsub}

\begin{proof} We will only prove {\bf 1}. The proof of {\bf 2} follows easily.

Let \m{x\in C}. Since \ \m{E_x\simeq r\ko_{Y,x}\oplus\ko_{C,x}}, we have \ 
\m{\Ext^1_{\ko_{Y,x}}(E_x,E_x)\simeq\Ext^1_{\ko_Y,x}(\ko_{C,x},\ko_{C,x})}.
We will give an explicit description of the extension
\[0\lra\ko_{C,x}\lra\kv\lra\ko_{C,x}\lra 0\]
corresponding to \
\m{\phi(\omega_{\ke^{[1]},\ke^{[2]}})(x)\in\Ext^1_{\ko_Y}(\ko_{C,x},\ko_{C,x})},
 and from \ref{self}, {\bf 1} will follow from the fact that \ 
\m{\kv\simeq\ko_{Y,x}} \ if \m{x\not\in\kz\cap C}, and \ \m{\kv\simeq\ki_{Y,1}} 
\ if \m{x\in\kz\cap C}.

Let \ \m{\tau_1,\tau_2\in\Ext^1_{\ko_{Y_2,x}}(E_x,E_x)} \ be 
the images of \m{\sigma_1,\sigma_2} respectively.
We have also an exact sequence
\[0\lra\Ext^1_{\ko_{Y,x}}(E_x,E_x)\lra\Ext^1_{\ko_{Y_2,x}}(E_x,E_x)\lra
\End(E_x) \ , \]
and \ \m{\tau_1-\tau_2\in\Ext^1_{\ko_{Y,x}}(E_x,E_x)}. In a neighbourhood of 
$x$ in $\kc$, \m{\ke^{[1]}} is isomorphic to \ \m{r\ko_\kc\oplus\ko_{\kc_1}}, 
and \m{\ke^{[2]}} is isomorphic to \ \m{r\ko_\kc\oplus\ko_{\kc_2}}. We can 
suppose that these isomorphisms are the same on $Y$. The exact sequence 
$(\ref{equ1})$ is the canonical exact sequence
\xmat{0\ar[r] & r\ko_{Y,x}\oplus\ko_{C,x}\ar[r]\fleq[d] & r\ko_{Y_2,x}\oplus
\ko_{\kc_i\cap Y_2,x}\ar[r]\fleq[d] & r\ko_{Y,x}\oplus\ko_{C,x}\ar[r]\fleq[d] &
0\\
& E_x & \ke^{[i]}_x & E_x}
Note that \ \m{\ko_{\kc_i\cap Y_2,x}=\ko_{\kc_i,x}/(\pi_1^2)}. We have 
\[\tau_1-\tau_2 \ \in \ \Ext^1_{\ko_{Y_2,x}}(\ko_{C,x},\ko_{C,x}) \ \subset 
\ \Ext^1_{\ko_{Y_2,x}}(E_x,E_x) \ . \]
We have \ \m{\tau_1-\tau_2=\eta_1-\eta_2}, where \ 
\m{\eta_i\in\Ext^1_{\ko_{Y_2,x}}(\ko_{C,x},\ko_{C,x})} \ is associated to the 
canonical exact sequence
\xmat{0\ar[r] & \ko_{C,x}\ar[r]\fleq[d] & \ko_{\kc_i,x}/(\pi_i^2)\ar[r] & 
\ko_{C,x}\ar[r]\fleq[d] & 0\\ & \pi_i.\ko_{\kc_i,x}/(\pi_i^2) & &
\ko_{\kc_i,x}/(\pi_i)}
If \m{\beta\in\ko_{\kc_i,x}}, the image of \m{\pi_i\beta} in 
\m{\ko_{\kc_i,x}/(\pi_i^2)} depends only on \m{\beta_{|C}}. So for every 
\m{\alpha\in\ko_{C,x}} we can define \m{\pi_i\alpha\in\ko_{\kc_i,x}/(\pi_i^2)}. 
Let
\[\kn \ = \ \big\{(\pi_1\alpha,-\pi_2\alpha)\in\ko_{\kc_i,x}/(\pi_1^2)
\times\ko_{\kc_2,x}/(\pi_2^2) \ ; \ \alpha\in\ko_{C,x}\big\} \ , \]
which is a sub-\m{\ko_{\kc,x}}-module of \ \m{\ko_{\kc_i,x}/(\pi_1^2)
\times\ko_{\kc_2,x}/(\pi_2^2)}. Let
\[\ku \ = \ \big[\ko_{\kc_i,x}/(\pi_1^2)\times\ko_{\kc_2,x}/(\pi_2^2)\big]/\kn 
\ . \]
The morphism
\[\xymatrix@R=5pt{\Phi:\ku\ar[r] & \ko_{C,x}\times\ko_{C,x}\\
(\alpha_1,\alpha_2)\fmaps[r] & (\alpha_{1|C},\alpha_{2|C})}\]
is surjective. We have
\[\ker(\Phi) \ = \ \big\{(\pi_1\lambda_1,\pi_2\lambda_2) \ ; \ 
\lambda_1,\lambda_2\in\ko_{C,x}\big\}/\kn \ . \]
We have \ \m{\ker(\Phi)\simeq\ko_{C,x}}, the isomorphism being defined by
\[\xymatrix@R=5pt{\nu:\ko_{C,x}\ar[r] & \ker(\Phi)\\ \alpha\fmaps[r] &
(\pi_1\alpha,0)=(0,\pi_2\alpha) .}\]
Hence we have an exact sequence
\begin{equation}\label{equ2}
0\lra\ko_{C,x}\lra\ku\lra\ko_{C,x}\oplus\ko_{C,x}\lra 0 \ .
\end{equation}
We have an inclusion
\[\xymatrix@R=5pt{\mu_1:\ko_{\kc_1,x}/(\pi_1^2)\ar[r] & \ku\\
\alpha_1\fmaps[r] & (\alpha_1,0) ,}\]
and similarly \ \m{\mu_2:\ko_{\kc_2,x}/(\pi_2^2)\hookrightarrow\ku}. We have a 
commutative diagram with exact rows
\xmat{0\ar[r] & \ko_{C,x}\ar[r]\fleq[d] & 
\ko_{\kc_1,x}/(\pi_1^2)\ar[r]\flinc[d]^{\mu_1} & \ko_{C,x}\ar[r]\flinc[d]^\mu & 
0\\ 0\ar[r] & \ko_{C,x}\ar[r]^-\nu & \ku\ar[r] & \ko_{C,x}\oplus\ko_{C,x}\ar[r] 
& 0 \ ,}
where $\mu$ is the inclusion in the first factor.

Let \ \m{\gamma\in\Ext^1_{Y_2,x}(\ko_{C,x}\oplus\ko_{C,x},\ko_{C,x})=
\Ext^1_{Y_2,x}(\ko_{C,x},\ko_{C,x})\oplus\Ext^1_{Y_2,x}(\ko_{C,x},\ko_{C,x})} \ 
associated to $(\ref{equ2})$. From the preceding diagram and prop. 4.3.1 of 
\cite{dr1b}, 
the first component of $\gamma$ is \m{\eta_1}. Similarly the second component 
of $\gamma$ is \m{\eta_2}. So we have \ \m{\gamma=(\eta_1,\eta_2)}. It follows 
that \m{\eta_1-\eta_2} corresponds to the top exact sequence in the following 
commutative diagram
\xmat{0\ar[r] & \ko_{C,x}\ar[r]\fleq[d] & \kv\ar[r]\flinc[d] & 
\ko_{C,x}\ar[r]\ar[d]^\psi & 0\\
0\ar[r] & \ko_{C,x}\ar[r] & \ku\ar[r] & \ko_{C,x}\oplus\ko_{C,x}\ar[r] & 0}
where \ \m{\psi:\alpha\mapsto(\alpha,-\alpha)}, and
\[\kv \ = \ \big\{(\alpha,\beta)\in\ku \ ; \ \alpha_{|C}+\beta_{|C}=0\big\} \ .
\]
If \ \m{(u_1,u_2)\in\ko_{\kc_i,x}/(\pi_1^2)\times\ko_{\kc_2,x}/(\pi_2^2)} \ is
such that \ \m{u_{1|C}+u_{2|C}=0}, we will denote by \m{[u_1,u_2]} the 
corresponding element of $\kv$.
 
We have \ \m{\pi\kv=\nsp}, i.e. the top exact sequence is a sequence of 
\m{\ko_{Y,x}}-modules.

\sepprop

{\em The case \ \m{x\not\in\kz\cap C} -- } We have then
\[\ko_{\kc,x} \ = \ \big\{(\alpha_1,\alpha_2)\in\ko_{\kc_1,x}\times\ko_{\kc_2,x}
\ ; \ \alpha_{1|C}=\alpha_{2|C}\big\} \ . \]
Let
\[\xymatrix@R=5pt{{\bf f}:\ko_{\kc,x}\ar[r] & \kv\\ 1\fmaps[r] & [1,-1] \ . }\]
We now prove that {\bf f} induces an isomorphism \ \m{\ko_{Y,x}\simeq\kv}.
It is obvious that {\bf f} is surjective and that \ \m{(\pi)\subset\ker({\bf 
f})}. Suppose that \m{(\alpha_1,\alpha_2)\in\ko_{\kc,x}} is such that \ \m{{\bf 
f}(\alpha_1,\alpha_2)=0}. We can then write
\[(\alpha_1,-\alpha_2) \ = \ 
(\pi_1\beta_1,-\pi_2\beta_2)+(\pi_1^2\epsilon_1,-\pi_2\epsilon_2^2) \ , \]
with \ \m{\beta_{1|C}=\beta_{2|C}}. Hence \ \m{(\beta_1,\beta_2)\in\ko_{\kc,x}},
and
\[(\alpha_1,\alpha_2) \ = \ 
\pi.\big[(\beta_1,\beta_2)+(\pi_1\epsilon_1,\pi_2\epsilon_2)\big] \ . \]
We have \ 
\m{(\beta_1,\beta_2)+(\pi_1\epsilon_1,\pi_2\epsilon_2)\in\ko_{\kc_x}}, 
hence \ \m{(\alpha_1,\alpha_2)\in(\pi)}.

\sepprop

{\em The case \ \m{x\in\kz\cap C} -- } We have then an isomorphism
\[\theta:\ko_{\kc_1,x}/(\pi_1\lambda_1)\lra\ko_{\kc_2,x}/(\pi_2\lambda_2)\]
such that \ \m{\theta(\pi_1)=\pi_2}, \m{\theta(\alpha)_{|C}=\alpha_{|C}} \ for 
every \ \m{\alpha\in\ko_{\kc_1,x}/(\pi_1\lambda_1)} (cf. \ref{prel2}). The 
restrictions \m{\lambda_{i|C}}, \m{i=1,2} are generators of the maximal ideal 
of \m{\ko_{C,x}}. We can also assume that \ \m{\theta(\lambda_1)=\lambda_2}. We 
have then
\begin{equation}\label{equ3}\ko_{\kc,x} \ = \ 
\big\{(\alpha_1,\alpha_2)\in\ko_{\kc_1,x}\times\ko_{\kc_2,x}
\ ; \ \theta(\alpha_1)=\alpha_2\big\} \ . \end{equation}
We now prove that \ \m{\kv\simeq\ki_x} (the ideal sheaf of \m{\{x\}} in 
$Y$). Let \ \m{z=(\pi_1\lambda_1,0)}, \m{t=(\lambda_1,\lambda_2)} in 
\m{\ko_{\kc,x}} (cf. $(\ref{equ3})$). We have \m{\ki_x=(t,z)}, $z$ is an 
equation of $C$ (in $Y$) and \m{t_{|C}} is a generator of the maximal ideal of 
\m{\ko_{C,x}}. Then there exists a unique morphism \ \m{\rho:\ki_x\to\kv} \ 
such that
\[\rho(t) \ = \ [1,-1] \ , \qquad \rho(z) \ = \ [\pi_1,0] \ . \]
To prove this we have only to show that if \m{\alpha,\beta\in\ko_{Y,x}} are 
such that \ \m{\alpha z+\beta t=0}, then we have \ 
\m{\beta[1,-1]+\alpha[\pi_1,0]=0} \ in $\kv$. We have \ \m{\alpha z+\beta t=0} 
\ if and only if we can write
\[\alpha \ = \ \epsilon t+\gamma z \ , \qquad \beta \ = \ -\epsilon z \ , \]
with \m{\epsilon,\gamma\in\ko_{Y,x}}. We have then
\begin{eqnarray*}
\beta[1,-1]+\alpha[\pi_1,0] & = & -\epsilon z[1,-1]+(\epsilon t+\gamma 
z)[\pi_1,0]\\
& = & \epsilon(\lambda_1,\lambda_2)[\pi_1,0]-\epsilon(\pi_1\lambda_1,0)[1,-1]+
\gamma[\pi_1^2\lambda_1,0]\\ & = & 0 \ .
\end{eqnarray*}
Now we show that $\rho$ is injective. Suppose that \ \m{(\alpha_1,\alpha_2),
(\beta_1,\beta_2)\in\ko_{\kc,x}} \ are such that
\[\rho((\alpha_1,\alpha_2)z+(\beta_1,\beta_2)t) \ = \ 0 \ . \]
Then we have
\[(\alpha_1,\alpha_2)[\pi_1,0]+(\beta_1,\beta_2)[1,-1] \ = \ 
[\alpha_1\pi_1+\beta_1,-\beta_2] \ = \ 0 \ . \]
Hence we can write
\[(\alpha_1\pi_1+\beta_1,-\beta_2) \ = \ (\pi_1\tau_1,-\pi_2\tau_2)+
(\pi_1^2\theta_1,\pi_2^2\theta_2) \ , \]
for some \m{\tau_i,\theta_i\in\ko_{\kc_i,x}} such that \ 
\m{\tau_{1|C}=\tau_{2|C}}, i.e.
\begin{equation}\label{equ4}
\alpha_1\pi_1+\beta_1 \ = \ \pi_1\tau_1+\pi_1^2\theta_1 \ , \qquad
\beta_2 \ = \pi_2\tau_2-\pi_2^2\theta_2 \ . 
\end{equation}
Let
\[u \ = \ (\alpha_1,\alpha_2)z+(\beta_1,\beta_2)t \ = \
\big(\lambda_1(\alpha_1\pi_1+\beta_1),\lambda_2\beta_2\big) \ . \]
From $(\ref{equ4})$ we see that \m{\beta_1} is a multiple of \m{\pi_1}, and 
\m{\beta_2} a multiple of \m{\pi_2}: \m{\beta_1=\pi_1\beta'_1}, 
\m{\beta_2=\pi_2\beta'_2}. We have then
\[u \ = \ \big(\pi_1\lambda_1(\alpha_1+\beta'_1),\pi_2\lambda_2\beta'_2\big) \ .
\]
We have 
\[\alpha_1+\beta'_1 \ = \ \tau_1+\pi_1\theta_1 \ , \
\beta'_2 \ = \ \tau_2-\pi_2\theta_2 \ . \]
Hence \ \m{(\alpha_1+\beta'_1)_{|C}=\beta'_{2|C}} \ and \ 
\m{\dsp\big(\lambda_1(\alpha_1+\beta'_1),\lambda_2\beta'_2\big)\in\ko_{\kc,x}}.
It follows that \m{u=0} in \m{\ko_{Y,x}}.

Now we show that $\rho$ is surjective. Let \ \m{[\alpha,\beta]\in\kv}. Then \ 
\m{\alpha_{|C}=-\beta_{|C}}. Let \m{\mu\in\ko_{\kc_1,x}} \ be such that \ 
\m{(\mu,-\beta)\in\ko_{\kc,x}}. We have
\[[\alpha,\beta]-(\mu,-\beta)[1,-1] \ = \ [\alpha-\mu,0] \ . \]
We can write \ \m{\alpha-\mu=\pi_1\zeta}, \m{\zeta\in\ko_{\kc_1,x}}. Let \ 
\m{\delta\in\ko_{\kc_2,x}} \ be such that \ \m{(\zeta,\delta)\in\ko_{\kc,x}}. 
We have then
\[[\alpha,\beta] \ = \ (\mu,-\beta)[1,-1]+(\zeta,\delta)[\pi_1,0] \ = \ 
\rho((\mu,-\beta)t+(\zeta,\delta)z) \ . \]
\end{proof}

\end{sub}

\vskip 1.5cm

\end{document}